# An Algebraic–Operator Construction of the Half-Derivative on a Graded Monomial Space

## Part I — Recurrence, Double Factorials, the Wallis Product, and Normalization

Davit Kapanadze
Georgian National University SEU, Tbilisi, Georgia
Associate Member, Institute of Mathematics and its Applications (IMA), United Kingdom
ORCID: 0009-0006-8112-5084
Corresponding author: d.kapanadze@seu.edu.ge

## Abstract

This paper constructs a half-order differentiation operator on the graded monomial space spanned by non-negative integer and half-integer powers of x, with $x > 0$. The algebraic stage uses neither an integral kernel, a limiting process, nor the Gamma function. The operator is assumed to act in the form

$$D^{1/2}x^{\beta} = c(\beta)x^{\beta-1/2},$$

Requiring two successive applications of the operator to reproduce the ordinary first derivative yields the fundamental recurrence relation

$$c(\beta)c(\beta-1/2) = \beta.$$

Expanding the recurrence produces two linked coefficient families, one for integer powers and one for half-integer powers. Their formulas contain a free normalization constant $c_0 = c(0)$, which cancels under composition. Thus, the compositional requirement fixes the relative coefficients but not the scale of each individual half-step. Ratios of even and odd double factorials lead naturally to a partial Wallis product; however, the Wallis product alone does not determine $c_0$. Imposing compatibility with the monomial formula for the left Riemann–Liouville half-derivative with lower limit 0 selects the value

$$c_0 = \frac{1}{\sqrt{\pi}}.$$

is selected. With this normalization,

$$D^{1/2} x^{\beta} = \frac{\Gamma(\beta+1)}{\Gamma\left(\beta+\frac{1}{2}\right)} x^{\beta-\frac{1}{2}}$$

holds for every $\beta \in \{0, 1/2, 1, 3/2, \ldots\}$. The resulting monomial formula agrees with the corresponding Riemann–Liouville formula. For positive integer powers it also agrees with the Caputo formula, whereas for the constant function it agrees only with the Riemann–Liouville rule. The resulting linear compositional operator lowers the exponent by exactly one half and is fully specified on the stated graded monomial space. Extension to broader function spaces, construction of an integral or convolution kernel, and the analysis of non-locality remain open problems.

## Introduction

The ordinary derivative is one of the fundamental concepts of mathematical analysis. Analytically, it describes the local rate of change of a function; geometrically, it gives the slope of the tangent line; and, from the operator viewpoint, it maps a function to its derivative. If the differentiation operator is denoted by D, then for a monomial

$$D x^n = n x^{n-1}.$$

Two effects occur simultaneously in this formula: the exponent is reduced by one, and the coefficient n appears. Repeated differentiation produces higher-order integer operators, and for suitable functions the composition law holds:

$$D^{\alpha} \circ D^{\beta} = D^{\alpha+\beta},$$

when α and β are non-negative integers.

The additivity of integer orders naturally suggests extending differentiation to non-integer orders. In particular, can one construct a half-order operator whose composition with itself equals the ordinary first derivative?

$$D^{1/2} \circ D^{1/2} = D \, ?$$

In this operator sense, $D^{1/2}$ may be regarded as a compositional square root of D. The term "square root" here does not mean a direct transfer of the numerical root operation to functions. It denotes an operator whose two successive applications on a given functional family yield the action of ordinary differentiation.

The extension of differentiation and integration to non-integer orders is the subject of fractional calculus. The question dates to the end of the seventeenth century, and its subsequent development is associated with Liouville, Riemann, Grünwald, Letnikov, Caputo, and many other authors [14, 18]. Among the most widely used classical formulations are the Riemann–Liouville, Caputo, and Grünwald–Letnikov operators [1–4, 20, 21]. Their definitions involve integral kernels, ordinary derivatives, finite-difference limits, the Gamma function, or combinations of these ingredients.

The fractional derivative is not a single universal operator. Different definitions may give different results for the same function. The result is influenced by the type of operator, its left- or right-sided direction, the lower limit of integration, the regularity of the function, the initial conditions, and the functional space under consideration. For example, when the lower limit is chosen to be 0, the Riemann–Liouville half-derivative of the constant function is nonzero:

$$D_{0+}^{RL,1/2} 1 = \frac{1}{\sqrt{\pi}} x^{-1/2},$$

whereas the Caputo half-derivative of the same constant equals zero:

$$D_{0+}^{C,1/2} 1 = 0.$$

The present paper approaches the half-order operator from a different direction. At the algebraic stage, no pre-existing integral definition, Gamma-function formula, or finite-difference limit is assumed. The starting point is the compositional requirement

$$\left(D^{1/2}\right)^2 = D.$$

Since the ordinary derivative reduces the power of a monomial by one, the natural starting form for the half-order operator is taken to be

$$D^{1/2} x^{\beta} = c(\beta)\, x^{\beta - \frac{1}{2}},$$

where c(β) is initially an unknown coefficient. Applying the operator a second time gives

$$D^{1/2}\left(D^{1/2}x^{\beta}\right)=c(\beta)\,c\left(\beta-\frac{1}{2}\right)x^{\beta-1}.$$

On the other hand,

$$D\,x^{\beta}=\beta\,x^{\beta-1}.$$

Comparing these expressions yields the fundamental recurrence relation of the construction:

$$c(\beta)\,c\left(\beta-\frac{1}{2}\right)=\beta$$

This recurrence is the algebraic foundation of the construction.

The investigation is restricted to linear operators that lower a monomial exponent by exactly one half and multiply each monomial by a scalar coefficient. Accordingly, the result is a compositional square root of differentiation within this specified operator class; it is not a classification of all possible square roots of D.

Expanding the recurrence produces two linked coefficient families for integer and half-integer powers and leaves one free normalization constant. Ratios of even and odd double factorials generate a partial Wallis product, but this product does not determine the normalization. A specific positive normalization is selected by requiring compatibility with the Riemann–Liouville monomial formula.

With this normalization, the coefficients agree with the corresponding Riemann–Liouville monomial coefficients; for positive integer powers they also agree with the Caputo formula. This monomial-level agreement does not imply equivalence of the full operator theories. The domain, the boundary term, and the limitations of the present construction are stated explicitly below.

The main tasks of the paper are:

1. formulation of the half-order compositional requirement on the graded family of monomials;

2. derivation of the recurrence relation for the coefficients;

3. construction of the coefficient families for integer and half-integer powers;

4. demonstration of the emergence of double factorials and the Wallis product;

5. separation of the compositional structure from the choice of normalization;

6. selection of a normalization compatible with the Riemann–Liouville formula;

7. verification of the obtained formulas on monomials and their finite linear combinations;

8. comparison of the results with the Riemann–Liouville and Caputo formulations;

9. clarification of the boundaries of the obtained construction and directions for further development.

The paper is organized as follows. Chapter I formulates the operator setting and the half-order composition requirement. Chapter II derives the recurrence and the integer and half-integer coefficient families. Chapter III develops the double-factorial and Wallis-product structure, selects the normalization, and verifies composition. Chapter IV compares the construction with classical fractional operators. Chapter V summarizes possible application domains of classical fractional calculus and states the additional results required before the present construction can be used in such models.

## Chapter I. The Derivative as an Operator and the Idea of the Half-Derivative

### 1.1. The Derivative as an Operator

Let us denote the differentiation operator by D:

$$D[f](x) = f'(x), \qquad Df = f'.$$

For a monomial the familiar rule holds:

$$Dx^n = nx^{n-1}.$$

The operator D is linear:

$$D(af + bg) = aDf + bDg.$$

Therefore, if

$$P(x) = \sum_{k=0}^{N} a_k x^k,$$

then

$$DP(x)=\sum_{k=1}^{N} k\,a_k\,x^{k-1}.$$

Below the same linearity principle will be used to extend the action of the half-order operator from monomials to their finite linear combinations.

### 1.2. Repeated Differentiation

Repeated application of the operator is defined by composition:

$$D^2 = D \circ D,$$

and, in general, for an m-fold composition,

$$D^m = D \circ D \circ \cdots \circ D \quad \text{(m times)}.$$

For a monomial, when $m \le n$,

$$D^m x^n = n(n-1)\cdots(n-m+1)\,x^{n-m} = \frac{n!}{(n-m)!}\,x^{n-m}.$$

If $m > n$, then after a sufficient number of differentiations the result becomes zero.

### 1.3. The Additivity Law for Operator Orders

The basic property of repeated differentiation is

$$D^m \circ D^r = D^{m+r},$$

where m and r are non-negative integers. Thus, the orders of differentiation add under composition. It is precisely this property that provides the initial motivation for extending the order of the operator from integer to non-integer values.

### 1.4. The Operator Idea of the Half-Derivative

Since

$$\frac{1}{2} + \frac{1}{2} = 1,$$

it is natural to seek an operator whose two successive applications yield the ordinary first derivative. We denote such an operator by $D^{1/2}$ and require

$$D^{1/2} \circ D^{1/2} = D.$$

Acting on a function, this means

$$D^{1/2}(D^{1/2}f) = Df.$$

In this sense, $D^{1/2}$ is an operator half-step toward one full differentiation. The composition condition alone, however, does not determine its action. We therefore begin with monomials, for which changes in exponent and coefficient can be tracked explicitly.

The analysis is restricted to linear operators that lower the monomial grading by one half-step. The conclusions therefore concern this class, not every possible operator square root of D.

### 1.5. Action of the Half-Order Operator on a Monomial

For the ordinary derivative it is known that

$$Dx^{\beta} = \beta x^{\beta-1}.$$

Assume that each application of the half-order operator lowers the exponent by 1/2:

$$D^{1/2}x^{\beta} = c(\beta)x^{\beta-1/2},$$

where $c(\beta)$ is an unknown coefficient. Applying the operator a second time gives

$$D^{1/2}(D^{1/2}x^{\beta}) = c(\beta)c(\beta-1/2)x^{\beta-1}.$$

From the compositional requirement follows the fundamental relation for the coefficients:

$$c(\beta)c(\beta-1/2) = \beta.$$

This recurrence underlies the entire construction. Before normalization, the notation $D^{1/2}$ represents a one-parameter family indexed by the nonzero initial coefficient $c_0 = c(0)$; Chapter III selects one member of this family.

Consider the family of monomials

$$1, x^{1/2}, x, x^{3/2}, \dots, x^{m/2}, \dots, \quad x > 0.$$

Let $V_+$ denote the finite real linear span of the monomials $1, x^{1/2}, x, x^{3/2}, \dots$ on $x > 0$. The restriction $x > 0$ ensures that all half-integer powers are single-valued and real. For the constant function,

$$D^{1/2}1 = c(0)x^{-1/2}.$$

Because the first half-step applied to 1 produces $x^{-1/2}$, we enlarge the initial space by adjoining this monomial. The extended space is

$$\tilde{V} = V_{+} \oplus \operatorname{span}\{x^{-1/2}\}.$$

The operator is defined as

$$D^{1/2} : \tilde{V} \to \tilde{V},$$

while the compositional requirement is imposed on the elements of the initial space:

$$(D^{1/2})^2 f = Df,$$

where f is an element of the initial space.

For the constant function the following must hold:

$$D^{1/2}(D^{1/2}1) = D1 = 0.$$

For $\beta = 0$, the same composition identity gives $c(0)c(-1/2) = 0$. Since $c(0) = c_0 \neq 0$, the boundary coefficient is forced to satisfy $c(-1/2) = 0$; equivalently,

$$D^{1/2}x^{-1/2} = 0.$$

Consequently,

$$D^{1/2}(D^{1/21}) = c_0 \cdot D^{1/2}x^{-1/2} = 0 = D1,$$

Thus, the extension to $\tilde{V}$ is consistent with the composition law on $V_{+}$, and no additional boundary parameter is introduced.

## 1.6. Methodological Position and Scope of the Construction

The present approach does not begin with a known integral formula. It starts from operator composition, derives the coefficient recurrence and the two graded coefficient families, introduces double factorials, and only then selects a normalization and compares the result with classical formulas. No limit, integral kernel, or Gamma function is used in the algebraic derivation.

The operator defined on monomials extends linearly to finite sums:

$$D^{1/2}\left(\sum_{j=1}^{N} a_j x^{\beta_j}\right) = \sum_{j=1}^{N} a_j c\left(\beta_j\right) x^{\beta_j - 1/2}.$$

The current objective is to determine the coefficient structure on a specific graded monomial space. General function spaces, integral or convolution kernels, a full theory of non-locality, and a semigroup law for arbitrary real orders are beyond the scope of this paper and require a separate functional extension.

# Chapter II. Algebraic Construction of the Half-Derivative

## 2.1. The Basic Principle of Composition

For ordinary differentiation,

$$Dx^{n} = nx^{n-1}.$$

For the half-order operator we take

$$D^{1/2}x^{\beta} = c(\beta)x^{\beta-1/2},$$

and require

$$D^{1/2}(D^{1/2}x^{\beta}) = Dx^{\beta}.$$

The coefficients for integer and half-integer powers must be determined together, because one half-step maps the integer sector to the half-integer sector. Their values are fixed, up to one normalization constant, by the composition requirement.

## 2.2. Derivation of the Recurrence

The first application is

$$D^{1/2}x^{\beta} = c(\beta)x^{\beta-1/2}.$$

The second application gives

$$D^{1/2}(D^{1/2}x^{\beta}) = c(\beta)c(\beta-1/2)x^{\beta-1}.$$

According to the compositional requirement, this must equal

$$Dx^{\beta} = \beta x^{\beta-1}.$$

Thus,

$$c(\beta)c(\beta-1/2) = \beta.$$

The recurrence follows solely from the operator composition requirement; no integral, limiting process, or Gamma function is used.

## 2.3. First Values of the Recurrence

Take as the initial value

$$c(0) = c_0,\ c_0 \neq 0.$$

At this stage $c_0$ is a free normalization constant.

If $\beta = 1/2$, then

$$c(1/2)c(0) = 1/2,\ \ c(1/2) = 1/(2c_0).$$

If $\beta = 1$, then

$$c(1)c(1/2) = 1,\ \ c(1) = 2c_0.$$

If $\beta = 3/2$, then

$$c(3/2)c(1) = 3/2,\ \ c(3/2) = 3/(4c_0).$$

If $\beta = 2$, then

$$c(2)c(3/2) = 2,\ \ c(2) = (8/3)c_0.$$

If $\beta = 5/2$, then

$$c(5/2)c(2) = 5/2,\ \ c(5/2) = 15/(16c_0).$$

If $\beta = 3$, then

$$c(3)c(5/2) = 3,\ \ c(3) = (16/5)c_0.$$

Thus, $c_0$ occurs in the integer-sector coefficients, whereas $1/c_0$ occurs in the half-integer-sector coefficients. These reciprocal factors cancel under two successive half-step actions.

### 2.4. The Pattern for Integer Arguments

For integer $n \geq 1$,

$$c(n)c(n-\frac{1}{2}) = n,$$

while

$$c(n-\frac{1}{2})c(n-1) = n-\frac{1}{2}.$$

Dividing these two equalities gives

$$\frac{c(n)}{c(n-1)} = \frac{n}{n-\frac{1}{2}} = \frac{2n}{2n-1}.$$

Therefore,

$$c(n) = \frac{2n}{2n-1} c(n-1).$$

Expanding the recurrence,

$$c(n) = \frac{2n}{2n-1} \cdot \frac{2n-2}{2n-3} \cdots \frac{2}{1} c_0,$$

that is,

$$c(n) = c_0 \cdot [2 \cdot 4 \cdot 6 \cdots 2n]/[1 \cdot 3 \cdot 5 \cdots (2n-1)].$$

Correspondingly,

$$D^{1/2} x^n = c_0 \frac{2 \cdot 4 \cdot 6 \cdots 2n}{1 \cdot 3 \cdot 5 \cdots (2n-1)} x^{n-1/2}.$$

## 2.5. The Appearance of Double Factorials

The even double factorial is defined as

$$(2n)!! = 2 \cdot 4 \cdot 6 \cdots 2n,$$

and the odd double factorial as

$$(2n-1)!! = 1 \cdot 3 \cdot 5 \cdots (2n-1).$$

With this notation,

$$c(n) = c_0 \cdot (2n)!!/(2n-1)!!.$$

Thus,

$$D^{1/2} x^n = c_0 \cdot [(2n)!!/(2n-1)!!] \cdot x^{n-1/2}.$$

For the case $n = 0$, the standard conventions

$$0!! = 1, \ (-1)!! = 1$$

are used, as a result of which the same formula gives $c(0) = c_0$.

## 2.6. Relation Between Double Factorials and Ordinary Factorials

For the even double factorial,

$$(2n)!! = 2^n n!.$$

Since

$$(2n)! = (2n-1)!!\cdot(2n)!!,$$

we have

$$(2n-1)!! = \frac{(2n)!}{2^n n!}.$$

Therefore

$$c(n) = c_0 \frac{2^n n!}{(2n)!/(2^n n!)} = c_0 \frac{4^n (n!)^2}{(2n)!}.$$

Thus,

$$c(n) = c_0 \frac{4^n (n!)^2}{(2n)!}.$$

Using the central binomial coefficient

$$\binom{2n}{n} = \frac{(2n)!}{(n!)^2}.$$

we may also write

$$c(n) = c_0 \frac{4^n}{\binom{2n}{n}}.$$

## 2.7. Verification on the Monomial x

Since

$$c(1) = 2c_0,$$

we have

$$D^{1/2}x = 2c_0 x^{1/2}.$$

Also,

$$c(\tfrac{1}{2}) = \frac{1}{2c_0},$$

so

$$D^{1/2}x^{1/2} = \frac{1}{2c_0}.$$

The second application gives

$$D^{1/2}(D^{1/2}x) = 2c_0 \cdot \frac{1}{2c_0} = 1 = Dx.$$

Thus, for x the composition requirement holds for every $c_0 \neq 0$.

## 2.8. Verification of the Construction on the Example of $x^2$

Since

$$c(2) = \frac{8}{3} c_0,$$

we have

$$D^{1/2} x^2 = \frac{8}{3} c_0 x^{3/2}.$$

Also,

$$c\left(\frac{3}{2}\right) = \frac{3}{4c_0},$$

so

$$D^{1/2}x^{3/2} = [3/(4c_0)]x.$$

The second application gives

$$D^{1/2}\left(D^{1/2} x^2\right) = \frac{8}{3} c_0 \cdot \frac{3}{4c_0} x = 2x = Dx^2.$$

## 2.9. General Monomial Formula and Summary of Chapter II

For a general positive integer n,

$$D^{1/2} x^n = c(n) x^{n-1/2}.$$

The second application gives

$$D^{1/2}\left(D^{1/2} x^n\right) = c(n)\, c\left(n - \frac{1}{2}\right) x^{n-1}.$$

By the recurrence,

$$c(n)c\left(n-\frac{1}{2}\right)=n,$$

So

$$D^{1/2}\left(D^{1/2}x^n\right)=nx^{n-1}=Dx^n$$

From the compositional condition

$$c(n)c\left(n-\frac{1}{2}\right)=n$$

it follows that

$$c\left(n-\frac{1}{2}\right)=\frac{n}{c(n)}.$$

Substituting the formula for the integer-power coefficient, we obtain

$$c\left(n-\frac{1}{2}\right)=\frac{n}{c_0\frac{(2n)!!}{(2n-1)!!}}=\frac{n(2n-1)!!}{c_0(2n)!!}.$$

Since

(2n)!! = (2n)(2n−2)!!,

Finally

$$c\left(n-\frac{1}{2}\right)=\frac{(2n-1)!!}{2c_0(2n-2)!!},\ \ n \geq 1.$$

For example, for n = 1:

$$c\left(\frac{1}{2}\right)=\frac{1!!}{2c_0 0!!}=\frac{1}{2c_0},$$

since 1!! = 1 and 0!! = 1.

Verification:

$$c(1)=c_0\frac{2!!}{1!!}=2c_0,$$

So

$$c(1)c\left(\frac{1}{2}\right)=2c_0\cdot\frac{1}{2c_0}=1,$$

which exactly matches the compositional requirement.

For the constant function we obtain

$$D^{1/2}1 = c_0 x^{-1/2}.$$

Using the boundary value obtained in Chapter I,

$$D^{1/2}x^{-1/2} = 0,$$

we get

$$D^{1/2}(D^{1/2}1) = 0 = D1.$$

By linear extension, if

$$P(x)=\sum_{k=0}^{N} a_k x^k,$$

then

$$D^{1/2}(D^{1/2}P) = DP.$$

Chapter II has derived the coefficient recurrence, identified the linked integer and half-integer families, and verified the composition law for every nonzero $c_0$. The normalization remains undetermined at this stage.

## Chapter III. The Normalization Issue, the Wallis Product, and Normalized Formulas for the Half-Derivative

The recurrence obtained in Chapter II,

$$c(\beta)c\left(\beta-\frac{1}{2}\right)=\beta$$

determines the mutual relation among the coefficients but does not fix their common scale. For integer powers we obtained

$$c(n)=c_0\frac{(2n)!!}{(2n-1)!!},$$

while for the corresponding half-integer powers

$$c\left(n-\frac{1}{2}\right)=\frac{(2n-1)!!}{2c_0(2n-2)!!},n\geq 1.$$

Chapter III explains how a Wallis-type product emerges from the even-to-odd double-factorial ratios, clarifies its connection with π, and separates that connection from the independent choice of normalization. The recurrence and coefficient families were obtained without limits, integrals, or the Gamma function. The classical Wallis limit is used only to identify the limiting numerical constant, and the Gamma function appears later only as a comparison notation.

## 3.1. Two Natural Ratios

Denote the product part of the integer-power coefficient by

$$A_n = \frac{(2n)!!}{(2n-1)!!}.$$

Then

$$c(n) = c_0 A_n.$$

In expanded form,

$$A_n = \frac{2 \cdot 4 \cdot 6 \cdots (2n)}{1 \cdot 3 \cdot 5 \cdots (2n-1)}.$$

Now introduce the related second ratio

$$B_n = \frac{(2n)!!}{(2n+1)!!}.$$

Since

$$(2n+1)!! = (2n+1)(2n-1)!!,$$

we have

$$B_n = \frac{A_n}{2n+1}.$$

The first ratio appears directly in the integer-sector coefficients. The second compares the same even product with an odd product containing one additional factor.

### 3.2. Emergence of the Wallis-Type Product

Consider

$$W_n = A_n B_n.$$

Then

$$W_n = \frac{\left((2n)!!\right)^2}{(2n-1)!!(2n+1)!!}$$

In expanded form,

$$W_n = \frac{2^2}{1\cdot 3}\cdot\frac{4^2}{3\cdot 5}\cdot\frac{6^2}{5\cdot 7}\cdots\frac{(2n)^2}{(2n-1)(2n+1)}$$

Thus,

$$W_n = \prod_{k=1}^{n}\frac{(2k)^2}{(2k-1)(2k+1)} = \prod_{k=1}^{n}\frac{4k^2}{4k^2-1}.$$

This is the nth partial product of the Wallis product. The normalization constant $c_0$ does not occur in $W_n$. Hence the Wallis product reveals a structural connection with π but cannot, by itself, determine $c_0$.

### 3.3. Numerical Behavior of $W_n$ and its Connection to π

By the classical Wallis formula [19],

$$\lim_{n\to\infty} W_n = \frac{\pi}{2}$$

equivalently,

$$\lim_{n\to\infty} 2W_n = \pi$$

The expanded infinite product is

$$\frac{\pi}{2} = \frac{2^2}{1\cdot 3}\cdot\frac{4^2}{3\cdot 5}\cdot\frac{6^2}{5\cdot 7}\cdots.$$

Here $W_n$ approaches π/2, whereas $2W_n$ approaches π. A proof of Wallis's formula is outside the scope of this paper. The relevant point is that the same product arises naturally from the even and odd factors in the half-derivative coefficients.

The logical sequence is

$$\text{recurrence} \to \text{double factorials} \to W_n \longrightarrow \frac{\pi}{2}.$$

However, the limiting Wallis product identifies only the value π/2; it does not determine the normalization constant.

### 3.4. Computational Illustration of the Wallis Product

To illustrate the numerical behavior, one can use a simple program that computes

$$W_n = \prod_{k=1}^{n} \frac{(2k)^2}{(2k-1)(2k+1)}$$

and then $2W_n$. Numerically, one starts with W = 1 and, for k = 1, 2, ..., n, multiplies W by the factor

$$\frac{(2k)^2}{(2k-1)(2k+1)} ;$$

finally compute 2W. The following are several numerical values obtained:

| n | $W_n$ | $2W_n$ |
|---|---|---|
| 1 | 1.3333333333 | 2.6666666667 |
| 2 | 1.4222222222 | 2.8444444444 |
| 5 | 1.5010879773 | 3.0021759546 |
| 10 | 1.5338519033 | 3.0677038066 |
| 100 | 1.5668937453 | 3.1337874906 |
| 1000 | 1.5704038730 | 3.1408077460 |
| 10000 | 1.5707570590 | 3.1415141180 |

This computation is illustrative rather than a proof; it shows numerically that $2W_n$ approaches π.

### 3.5. Choice of the Normalization Constant

The recurrence condition obtained in Chapter II,

$$c(\beta)c(\beta-1/2) = \beta,$$

determines the dependence between neighboring half-step coefficients but leaves their common normalization free. For integer powers,

$$c(n)=c_0\frac{(2n)!!}{(2n-1)!!},$$

while for the corresponding half-integer powers

$$c\left(n-\frac{1}{2}\right)=\frac{(2n-1)!!}{2c_0(2n-2)!!},n\geq 1$$

These formulas satisfy the composition requirement for every $c_0 \neq 0$. In particular,

$$c(n)c\left(n-\frac{1}{2}\right)=n.$$

The normalization cancels because

$$c_0\cdot\frac{1}{c_0}=1.$$

Therefore, $c_0$ does not affect the result of two successive half-order actions:

$$D^{1/2} \circ D^{1/2} = D$$

on the initial graded space. It does, however, change each individual half-step. Replacing $c_0$ by $\lambda c_0$, with $\lambda \neq 0$, multiplies the integer-sector coefficients by $\lambda$ and the corresponding half-integer-sector coefficients by $\lambda^{-1}$. Different normalizations therefore produce different intermediate representations but the same final first derivative.

To select a specific operator, an additional normalization condition is required. We choose the positive normalization that agrees with the standard monomial formula for the left Riemann–Liouville half-derivative with lower limit 0. For the constant function, the classical coefficient is $1/\sqrt{\pi}$, as reviewed in Chapter IV; hence we set

$$c_0=\frac{1}{\sqrt{\pi}}.$$

The positive sign preserves positivity of the coefficients on positive monomials for x > 0. With this normalization,

$$D^{1/2}1=\frac{1}{\sqrt{\pi}}x^{-1/2}.$$

The value $c_0 = 1/\sqrt{\pi}$ is not implied by the recurrence or by the Wallis product. It is selected by the additional compatibility condition. Wallis's formula nevertheless provides the convergent numerical approximation

$$c_{0,n}=\frac{1}{\sqrt{2W_n}},$$

then

$$c_{0,n}\longrightarrow\frac{1}{\sqrt{\pi}}\approx 0.5641895835.$$

Thus, the recurrence determines the coefficient structure, whereas $c_0$ determines the relative scale of an individual half-step. The reciprocal normalization factors cancel under composition, and $c_0 = 1/\sqrt{\pi}$ is selected by compatibility with the classical monomial normalization.

### 3.6. Normalized Coefficients for Integer Powers

After the chosen normalization,

$$c(n)=\frac{(2n)!!}{(2n-1)!!\sqrt{\pi}},n=0,1,2,\ldots,$$

Moreover,

$$D^{1/2}x^n=\frac{(2n)!!}{(2n-1)!!\sqrt{\pi}}x^{n-\frac{1}{2}}.$$

The first values are

$$D^{1/2}1=\frac{1}{\sqrt{\pi}}x^{-1/2},$$

$$D^{1/2}x=\frac{2}{\sqrt{\pi}}x^{1/2},$$

$$D^{1/2}x^2=\frac{8}{3\sqrt{\pi}}x^{3/2},$$

$$D^{1/2}x^3=\frac{16}{5\sqrt{\pi}}x^{5/2}.$$

In terms of ordinary factorials,

$$c(n)=\frac{4^n(n!)^2}{(2n)!\sqrt{\pi}}.$$

If

$$P(x)=\sum_{k=0}^{N} a_k x^k,$$

then by linear extension

$$D^{1/2}P(x)=\sum_{k=0}^{N} a_k \frac{(2k)!!}{(2k-1)!!\sqrt{\pi}} x^{k-\frac{1}{2}}.$$

### 3.7. Normalized Coefficients for Half-Integer Powers

From the compositional requirement

$$c(n)c\left(n-\frac{1}{2}\right)=n$$

we obtain

$$c\left(n-\frac{1}{2}\right)=\frac{n}{c(n)}.$$

Substituting the normalized integer-power coefficient,

$$c\left(n-\frac{1}{2}\right)=\frac{(2n-1)!!\sqrt{\pi}}{2(2n-2)!!},n=1,2,3,\ldots.$$

Therefore

$$D^{1/2}x^{n-\frac{1}{2}}=\frac{(2n-1)!!\sqrt{\pi}}{2(2n-2)!!}x^{n-1}.$$

Setting k = n − 1, the same formula takes the form

$$D^{1/2}x^{k+\frac{1}{2}}=\frac{(2k+1)!!\sqrt{\pi}}{2(2k)!!}x^k,k=0,1,2,\ldots.$$

The first values are

$$D^{1/2}x^{1/2}=\frac{\sqrt{\pi}}{2},$$
$$D^{1/2}x^{3/2}=\frac{3\sqrt{\pi}}{4}x,$$
$$D^{1/2}x^{5/2}=\frac{15\sqrt{\pi}}{16}x^2,$$
$$D^{1/2}x^{7/2}=\frac{35\sqrt{\pi}}{32}x^3.$$

For comparison with classical notation, the same coefficient can be rewritten using the Gamma function. This rewriting is not part of the derivation, and extension to general real powers is deferred to Part II:

$$c\left(n-\frac{1}{2}\right)=\frac{\Gamma\left(n+\frac{1}{2}\right)}{\Gamma(n)}.$$

The Gamma-function expression does not enter the derivation of the recurrence, the double-factorial formulas, or the normalized coefficients; it is only an equivalent classical notation for the result already obtained.

### 3.8. Table of Formulas: Obtaining the Ordinary Derivative from Two Half-Derivatives

For monomials of integer powers:

| Initial monomial | First half-derivative | Second half-derivative | Ordinary derivative |
|---|---|---|---|
| 1 | $\frac{1}{\sqrt{\pi}}x^{\frac{-1}{2}}$ | 0 | 0 |
| x | $\frac{2}{\sqrt{\pi}}x^{\frac{1}{2}}$ | 1 | 1 |
| $x^2$ | $\frac{8}{3\sqrt{\pi}}x^{\frac{3}{2}}$ | 2x | 2x |
| $x^3$ | $\frac{16}{5\sqrt{\pi}}x^{\frac{5}{2}}$ | $3x^2$ | $3x^2$ |

For monomials of half-integer powers:

| Initial monomial | First half-derivative | Second half-derivative | Ordinary derivative |
|---|---|---|---|
| $x^{\frac{1}{2}}$ | $\frac{\sqrt{\pi}}{2}$ | $\frac{1}{2}x^{\frac{-1}{2}}$ | $\frac{1}{2}x^{\frac{-1}{2}}$ |
| $x^{\frac{3}{2}}$ | $\frac{3\sqrt{\pi}}{4}x$ | $\frac{3}{2}x^{\frac{1}{2}}$ | $\frac{3}{2}x^{\frac{1}{2}}$ |
| $x^{\frac{5}{2}}$ | $\frac{15\sqrt{\pi}}{16}x^2$ | $\frac{5}{2}x^{\frac{3}{2}}$ | $\frac{5}{2}x^{\frac{3}{2}}$ |
| $x^{\frac{7}{2}}$ | $\frac{35\sqrt{\pi}}{32}x^3$ | $\frac{7}{2}x^{\frac{5}{2}}$ | $\frac{7}{2}x^{\frac{5}{2}}$ |

In both tables one can verify directly that

$$D^{1/2}\left(D^{1/2}x^{\beta}\right)=\beta x^{\beta-1}=D x^{\beta}.$$

### 3.9. The Half-Derivative of a Polynomial and the Compositional Property

Consider

$$P(x)=\sum_{k=0}^{N} a_k x^k.$$

By linearity,

$$D^{1/2}P(x)=\sum_{k=0}^{N} a_k \frac{(2k)!!}{(2k-1)!!\sqrt{\pi}} x^{k-\frac{1}{2}}.$$

After the second application—using $D^{1/2}x^{-1/2} = 0$ for the term $k = 0$ and the half-integer coefficient formula for $k \geq 1$—we obtain

$$D^{1/2}\left(D^{1/2}P(x)\right)=\sum_{k=1}^{N} k a_k x^{k-1}=DP(x).$$

For example, if

$$P(x) = a_0 + a_1x + a_2x^2 + a_3x^3,$$

then

$$D^{1/2}P(x)=\frac{a_0}{\sqrt{\pi}}x^{-1/2}+\frac{2a_1}{\sqrt{\pi}}x^{1/2}+\frac{8a_2}{3\sqrt{\pi}}x^{3/2}+\frac{16a_3}{5\sqrt{\pi}}x^{5/2}.$$

The second half-step application gives

$$D^{1/2}\left(D^{1/2}P(x)\right)=a_1+2a_2x+3a_3x^2=DP(x).$$

### 3.10. Summary of Chapter III

Chapter III has put the coefficient system into normalized form. The partial Wallis product displays a structural connection with π but does not determine the free normalization constant. Compatibility with the Riemann–Liouville monomial formula selects the positive normalization $c_0 = 1/\sqrt{\pi}$. The tables and polynomial example then verify the composition law on the initial graded space. No limit, integral, or Gamma function is used in deriving the recurrence or the double-factorial formulas.

## Chapter IV. Historical–Theoretical Context of Fractional Calculus and Positioning of the Obtained Construction

### 4.1. Purpose and Scope of the Chapter

The preceding chapters derived the coefficient recurrence, the double-factorial structure, and a normalization compatible with the Riemann–Liouville monomial formula. This chapter places the algebraic–operator construction in its classical context and addresses the following questions:

1. which classical monomial formula the obtained result coincides with;

2. what the difference is between the Riemann–Liouville and Caputo half-derivatives;

3. why, in the case of $x^2$, both classical formulations coincide with the presented construction;

4. how the compositional motivation relates to other fractional operators;

5. what the precise theoretical status of the obtained result is.

Classical integral and limiting formulas are introduced here only for comparison. They are not used in the derivation of Chapters I–III. Agreement on monomials does not imply equality of kernels, initial conditions, non-locality properties, or domains; monomial agreement must therefore be distinguished from equivalence of operator theories.

## 4.2. The Classical Idea of the Fractional Derivative and its Historical Development

Fractional calculus extends differentiation and integration from integer to non-integer orders. Its development, beginning with the late-seventeenth-century question and continuing through the work of Liouville, Riemann, Grünwald, Letnikov, Caputo, and others, is reviewed in [14, 18]. Classical formulations use integral kernels, finite-difference limits, transforms, and special functions [1–4].

Because definitions differ in direction, lower limit, regularity assumptions, initial conditions, and function space, fractional differentiation is best viewed as a family of related operators rather than a single universal formula.

## 4.3. The Riemann–Liouville Derivative and the Monomial Formula

Let $0 < \alpha < 1$. For a suitable function, the left Riemann–Liouville derivative with lower limit 0 is defined by the formula [3, 4, 10]

$$D_{0+}^{\mathrm{RL},\alpha} f(x) = \frac{d}{dx}\left[\frac{1}{\Gamma(1-\alpha)}\int_0^x (x-t)^{-\alpha} f(t)\,dt\right]$$

For a monomial with $\beta > -1$ and $x > 0$, the standard formula is

$$D_{0+}^{\mathrm{RL},\alpha} x^{\beta} = \frac{\Gamma(\beta+1)}{\Gamma(\beta+1-\alpha)} x^{\beta-\alpha}$$

For $\alpha = 1/2$,

$$D_{0+}^{\mathrm{RL},\frac{1}{2}} x^{\beta} = \frac{\Gamma(\beta+1)}{\Gamma\left(\beta+\frac{1}{2}\right)} x^{\beta-\frac{1}{2}}$$

If n is a non-negative integer, then

$$D_{0+}^{\mathrm{RL},\frac{1}{2}} x^{n} = \frac{\Gamma(n+1)}{\Gamma\left(n+\frac{1}{2}\right)} x^{n-\frac{1}{2}}$$

This coincides exactly with the normalized formula obtained in Chapter III. For the constant function,

$$D_{0+}^{\mathrm{RL},\frac{1}{2}} 1 = \frac{1}{\sqrt{\pi}} x^{-\frac{1}{2}}$$

It was precisely on the basis of compatibility with this classical value that Chapter III chose

$$c_0 = \frac{1}{\sqrt{\pi}}$$

The Riemann–Liouville formula is not used to derive the recurrence. It enters only to select the free normalization and to compare the resulting monomial coefficients with a classical operator.

### 4.4. The Caputo Derivative and the Case of the Constant Function

Let $0 < \alpha < 1$ and let f have suitable regularity, for example $f \in AC[0, b]$. The left Caputo derivative with lower limit 0 is defined by [4, 10, 20]

$$D_{0+}^{\mathrm{C},\alpha} f(x) = \frac{1}{\Gamma(1-\alpha)} \int_0^x (x-t)^{-\alpha} f'(t)\,dt$$

For the constant function,

$$D_{0+}^{\mathrm{C},\alpha} 1 = 0$$

while for Riemann–Liouville,

$$D_{0+}^{\mathrm{RL},\alpha} 1 = \frac{x^{-\alpha}}{\Gamma(1-\alpha)}$$

Thus, on the constant function, the present algebraic–operator construction agrees with the Riemann–Liouville rule and not with the Caputo rule.

For positive integer powers, $n \geq 1$, both classical formulations give the same result:

$$D_{0+}^{\mathrm{RL},\alpha} x^n = D_{0+}^{\mathrm{C},\alpha} x^n = \frac{\Gamma(n+1)}{\Gamma(n+1-\alpha)} x^{n-\alpha}$$

In general, for $0 < \alpha < 1$, the relation between the Riemann–Liouville and Caputo derivatives is

$$D_{0+}^{\mathrm{RL},\alpha} f(x) = D_{0+}^{\mathrm{C},\alpha} f(x) + \frac{f(0)}{\Gamma(1-\alpha)} x^{-\alpha}$$

Consequently, the two derivatives agree whenever $f(0) = 0$. In particular, this condition holds for $x^n$ with $n \geq 1$.

### 4.5. Comparative Example: The Half-Derivative of $x^2$ by Three Methods

Consider

$$f(x) = x^2, \quad x > 0.$$

### 4.5.1. By the Riemann–Liouville Method

By the definition of the Riemann–Liouville half-derivative,

$$D_{0+}^{\mathrm{RL},\frac{1}{2}}x^2=\frac{1}{\sqrt{\pi}}\frac{d}{dx}\int_0^x\frac{t^2}{\sqrt{x-t}}\,dt$$

With the substitution t = xu, dt = x du,

$$\int_0^x\frac{t^2}{\sqrt{x-t}}\,dt=x^{\frac{5}{2}}\int_0^1u^2(1-u)^{-\frac{1}{2}}\,du$$

The last integral equals

$$\int_0^1u^2(1-u)^{-\frac{1}{2}}\,du=\frac{16}{15}$$

Therefore

$$\int_0^x\frac{t^2}{\sqrt{x-t}}\,dt=\frac{16}{15}x^{\frac{5}{2}}$$

Differentiating, we obtain

$$D_{0+}^{\mathrm{RL},\frac{1}{2}}x^2=\frac{1}{\sqrt{\pi}}\frac{d}{dx}\left(\frac{16}{15}x^{\frac{5}{2}}\right)=\frac{8}{3\sqrt{\pi}}x^{\frac{3}{2}}$$

### 4.5.2. By the Caputo Method

Since

$$f'(t) = 2t,$$

we have

$$D_{0+}^{\mathrm{C},\frac{1}{2}}x^2=\frac{1}{\sqrt{\pi}}\int_0^x\frac{2t}{\sqrt{x-t}}\,dt$$

With the substitution t = xu,

$$\int_0^x\frac{2t}{\sqrt{x-t}}\,dt=2x^{\frac{3}{2}}\int_0^1u(1-u)^{-\frac{1}{2}}\,du$$

Here

$$\int_0^1u(1-u)^{-\frac{1}{2}}\,du=\frac{4}{3}$$

Therefore

$$D_{0+}^{\mathrm{C},\frac{1}{2}} x^2 = \frac{8}{3\sqrt{\pi}} x^{\frac{3}{2}}$$

### 4.5.3. By the Algebraic–Operator Construction

By the formula obtained in Chapter III,

$$D_{\mathrm{alg}}^{\frac{1}{2}} x^n = \frac{(2n)!!}{(2n-1)!!\sqrt{\pi}} x^{n-\frac{1}{2}}$$

For n = 2,

$$D_{\mathrm{alg}}^{\frac{1}{2}} x^2 = \frac{4!!}{3!!\sqrt{\pi}} x^{\frac{3}{2}} = \frac{8}{3\sqrt{\pi}} x^{\frac{3}{2}}$$

For the second half-step application,

$$D_{\mathrm{alg}}^{\frac{1}{2}} x^{\frac{3}{2}} = \frac{3\sqrt{\pi}}{4} x$$

Therefore

$$D_{\mathrm{alg}}^{\frac{1}{2}}\left(D_{\mathrm{alg}}^{\frac{1}{2}} x^2\right) = \frac{8}{3\sqrt{\pi}} \cdot \frac{3\sqrt{\pi}}{4} x = 2x = Dx^2$$

Accordingly,

$$D_{0+}^{\mathrm{RL},\frac{1}{2}} x^2 = D_{0+}^{\mathrm{C},\frac{1}{2}} x^2 = D_{\mathrm{alg}}^{\frac{1}{2}} x^2 = \frac{8}{3\sqrt{\pi}} x^{\frac{3}{2}}$$

This agreement is consistent with $x^2(0) = 0$. It does not imply that the three operators are identical on all functions; the constant function already distinguishes the Riemann–Liouville and Caputo rules.

## 4.6. The Grünwald–Letnikov Formulation and the Law of Composition

The typical form of the left Grünwald–Letnikov derivative is

$$D_{0+}^{\mathrm{GL},\alpha} f(x) = \lim_{h\to 0^+} h^{-\alpha} \sum_{k=0}^{\lfloor x/h \rfloor} (-1)^k \binom{\alpha}{k} f(x-kh)$$

when the limit exists and the corresponding convergence conditions are satisfied. Under suitable conditions this formula is connected to the Riemann–Liouville derivative [21].

Classical theory first defines a complete fractional operator and then studies its properties. Here the starting requirement is imposed only on a selected graded monomial space, from which the coefficient recurrence is derived. The result is therefore not a substitute for the full Grünwald–Letnikov, Riemann–Liouville, or Caputo theories. What is established is the composition identity

$$\left(D_{\mathrm{alg}}^{\frac{1}{2}}\right)^2 f = Df, \qquad f \in V_+$$

on the graded monomial space considered here. This identity, together with the coefficient recurrence, is the precise scope of the present result.

## 4.7. Central and Parametric Fractional Operators

Central fractional differences and derivatives provide another composition-related direction [22]. They are defined through bidirectional action, Fourier symbols, and convolution kernels. The present construction instead starts from half-step lowering of monomial exponents and a coefficient recurrence. Central operators are therefore a methodological parallel, not a direct equivalent.

Tempered fractional calculus provides an example in which an additional parameter modifies kernel decay, memory duration, and frequency response [23]. The constant $c_0$ has a different role: it rescales the two graded sectors reciprocally and cancels after two half-step actions.

A change in $c_0$ should therefore not be identified with a tempering parameter. On the present graded space it is an algebraic renormalization; any independent physical interpretation would require a broader functional model containing additional structure.

## 4.8. Theoretical Status of the Obtained Result

Established:

• the recurrence for the coefficients on the graded family of monomials;

• the existence of a free nonzero normalization constant;

• two mutually linked coefficient families for integer and half-integer powers;

• obtaining the ordinary derivative, on the initial space $V_+$, for monomials and their finite linear combinations, by composing two half-step actions;

• agreement with the Riemann–Liouville monomial formula when $c_0 = 1/\sqrt{\pi}$ is chosen;

• agreement with the Caputo monomial formula for positive integer powers;

• correspondence with the Riemann–Liouville rule, and difference from the Caputo rule, with respect to the constant function.

Not yet established:

• definition of the operator on a general functional space;

• existence of an integral or convolution kernel;

• general properties of causality and non-locality;

• complete Laplace and Fourier transform formulas;

• the semigroup law for arbitrary real order;

• the natural form of initial and boundary conditions;

• an independent physical meaning of the normalization constant in applied models.

Accordingly, the result is a compositional operator in a specified class of linear half-step power-lowering operators on a graded monomial space, extended by linearity to finite sums.

### 4.9. Methodological Significance and Chapter Summary

The methodological contribution is that the half-order coefficient structure is derived from ordinary differentiation, monomial grading, composition, recurrence, and double factorials. Integral and Gamma-function formulas enter only at the normalization and comparison stages. The construction is related to the author's earlier algebraic–geometric work on ordinary differentiation [15–17].

Comparison with classical theory shows agreement with the Riemann–Liouville monomial formula for the considered non-negative integer powers and with the Caputo formula for positive integer powers. This identifies the position of the construction without producing a general non-local operator theory.

## Chapter V. Applied Significance of Fractional Calculus and Perspectives for Development of the Present Construction

### 5.1. Purpose and Scope of the Applied Overview

Fractional-order operators are used to model systems whose current state depends on more than the current value alone. Classical fractional operators combine past values through integral kernels, weighted differences, or transform symbols and are therefore used to represent memory, hereditary effects, delayed response, anomalous transport, and multiscale dynamics [5–13, 24].

The applications reviewed here concern established Riemann–Liouville, Caputo, Grünwald–Letnikov, Riesz, and related operators. They are not yet applications of the monomial construction developed in this paper. The purpose of the chapter is to identify possible directions and the theoretical or numerical results required before such applications can be considered.

### 5.2. Memory, Non-Locality, and Multiscale Response

The ordinary first derivative is local: it describes instantaneous change at a point. Many physical and biological systems, however, depend on their history. Material deformation may depend on earlier loads, particle displacement on the preceding trajectory, biological response on accumulated stimuli, and drug concentration on the history of absorption and elimination.

In Riemann–Liouville and Caputo operators, past influence is represented by an integral kernel. A power-law kernel assigns different weights to different past times, while the non-integer order interpolates between classical dynamical regimes.

This memory interpretation does not automatically apply to the present monomial construction. Until an integral or convolution kernel, transform symbol, or another non-local representation is established, the operator should be viewed only as a compositional operator on a graded monomial space.

### 5.3. Viscoelasticity and the Half-Order Oscillator

Viscoelasticity is one of the best-known application areas of fractional derivatives [5, 7]. A viscoelastic material exhibits both elastic and dissipative behavior, and its response may depend on the loading history.

A simple form of the classical fractional constitutive model is

$$\sigma(t)=\eta_\alpha D_t^\alpha \varepsilon(t), 0<\alpha<1$$

where $\sigma(t)$ is the stress, $\varepsilon(t)$ is the strain, and $\eta_\alpha$ is a corresponding material parameter.

The typical equation of a single-degree-of-freedom fractional oscillator is [5, 24]

$$mx''(t) + cD_t^\alpha x(t) + kx(t) = F(t).$$

For the half-order it takes the form

$$m x''(t) + c D_t^{\frac{1}{2}} x(t) + kx(t) = F(t)$$

This oscillator is a possible test problem only after the present operator has been defined on an appropriate space of time-dependent functions and its transform, kernel, causality, and initial-condition structure have been established.

### 5.4. Anomalous Diffusion and Transport

Classical diffusion is often described by the equation

$$\frac{\partial u}{\partial t} = K \frac{\partial^2 u}{\partial x^2}$$

where the mean square displacement is proportional to time:

$$\langle x^2(t) \rangle \propto t$$

In a complex medium one may instead have

$$\langle x^2(t) \rangle \propto t^\alpha.$$

If $0 < \alpha < 1$, the process is termed subdiffusion; values $\alpha > 1$ correspond to superdiffusive scaling. A typical time-fractional diffusion model is

$$D_t^\alpha u(x,t) = K \frac{\partial^2 u}{\partial x^2}, 0 < \alpha < 1$$

Such models describe particle trapping, long waiting times, and memory-dependent transport [6]. Spatially fractional advection–dispersion equations model transport in heterogeneous media [13].

The possible connection of the present construction to this direction requires a considerably broader theory, since in diffusion models the essential elements are mass conservation, positivity, existence of a fundamental solution, initial conditions, and transform properties.

### 5.5. Biological and Pharmacokinetic Models

Fractional models are used in biology because many biological responses display multiscale temporal behavior, delays, and hereditary effects [8, 11]. Applications include tissue mechanics, biological-signal dynamics, cellular transport, pulmonary and cardiovascular models, and drug absorption, distribution, and elimination.

The classical first-order elimination model

$$\frac{dC}{dt} = -kC$$

can be replaced by the fractional form

$$D_t^{\alpha} C(t) = -kC(t), 0 < \alpha < 1$$

Solutions of such models often involve the Mittag–Leffler function rather than the exponential function, producing slower decay and a longer memory tail [9].

There is currently no basis for interpreting the normalization freedom of the present construction as an independent tissue, drug, or biological parameter. Such an interpretation would require a separately formulated and testable extended model.

### 5.6. Fractional Linear Systems, Signals, and Control

In fractional linear systems, non-integer-order operators generate factors of the form $s^{\alpha}$ in transfer functions [12, 24]. For example,

$$G(s) = \frac{1}{a_0 + a_1 s^{\alpha_1} + \cdots + a_n s^{\alpha_n}}$$

Such a representation allows the study of impulse and step responses, frequency characteristics, phase shift, and stability.

The classical form of the fractional-order PIλDμ controller is [25]

$$C(s) = K_p + K_i s^{-\lambda} + K_d s^{\mu}$$

The exponents λ and μ need not be integers, which can provide additional flexibility in shaping the frequency response [12, 24, 25].

Before applying the present construction in such systems, it must first be established whether a corresponding transform formula can be obtained for it, for example

$$L\{D^{\frac{1}{2}}f\}(s)=s^{\frac{1}{2}}F(s)-I_0(s)$$

Here $I_0(s)$ denotes terms associated with initial data, or a modified analogue. This expression illustrates only a possible transform structure; it is not an established formula for the present operator.

## 5.7. Normalization Freedom as Algebraic Similarity

For comparison with the classical Riemann–Liouville normalization, let us introduce the dimensionless ratio

$$\rho=c_0\sqrt{\pi}$$

Then

$$c_0=\frac{\rho}{\sqrt{\pi}}$$

and the classical choice corresponds to ρ = 1.

For integer powers one may write

$$D_\rho^{\frac{1}{2}}x^n=\rho D_1^{\frac{1}{2}}x^n$$

while on the corresponding half-integer powers

$$D_\rho^{\frac{1}{2}}x^{n-\frac{1}{2}}=\rho^{-1}D_1^{\frac{1}{2}}x^{n-\frac{1}{2}}$$

Therefore

$$D_\rho^{\frac{1}{2}}\left(D_\rho^{\frac{1}{2}}x^n\right)=Dx^n$$

On the first half-step, ρ appears; on the second, $\rho^{-1}$; under composition they cancel one another.

This family can be given a more precise interpretation at the algebraic level. Define the sector renormalization operator as follows:

$$S_\rho x^n = \rho x^n, S_\rho x^{n+\frac{1}{2}} = x^{n+\frac{1}{2}}$$

Then

$$D_\rho^{\frac{1}{2}} = S_\rho^{-1} D_1^{\frac{1}{2}} S_\rho$$

Operators with different values of ρ are therefore related by a similarity transformation and represent reciprocal rescalings of the two graded sectors. On the present space, ρ has no independent algebraic observable content. A physical interpretation could arise only in an extended model containing additional structure that is not preserved by this similarity.

Thus, within the current construction, different normalizations describe algebraically equivalent intermediate representations of the same first-order action.

### 5.8. Symbolic Verification of the Normalization Family

Before seeking a full kernel representation, one may verify the normalization family symbolically and numerically on finite half-integer-power sums:

$$f(x) = \sum_{m=0}^{M} a_m x^{\frac{m}{2}}$$

For various $\rho \neq 0$, one computes

$$g_\rho(x) = D_\rho^{\frac{1}{2}} f(x)$$

then

$$h_\rho(x) = D_\rho^{\frac{1}{2}} g_\rho(x)$$

and one checks whether

$$h_\rho(x) = Df(x)$$

for the terms on which the composition is defined.

One may compare the first-step coefficients, the intermediate graphs, the second-step result, and numerical error. Such an experiment does not establish memory or a physical effect; it only illustrates that similarity-related normalizations lead to the same final derivative through differently scaled intermediate representations.

### 5.9. Requirements to be Verified Before Application

Before the present construction can be used in an applied model, the following issues must be resolved:

- the precise domain of definition of the operator and the functional space;
- properties of linearity, closedness, and continuity;
- existence of an integral or convolution kernel;
- the nature of causality or bidirectional action;
- formulas for the Laplace and Fourier transforms;
- the conditions under which the index law holds;
- the natural form of initial and boundary conditions;
- existence of an inverse operator;
- solvability and uniqueness of the corresponding differential equations;
- stability and convergence of the numerical scheme;
- dimensional consistency of the normalization parameter;
- its identifiability with respect to other model parameters.

The theory of fractional differential equations shows that formal operator notation alone does not ensure well-posedness [10]. Work on central and tempered operators likewise highlights the importance of composition, convergence, and transform properties [21–23].

### 5.10. Summary of Applied Significance

These examples illustrate the established applications of fractional calculus: hereditary response in viscoelasticity [5, 7], non-classical scaling in anomalous diffusion and transport [6, 13], multiscale and delayed processes in biology and pharmacokinetics [8, 9, 11], and flexible frequency response in signals and control [12, 24, 25]. The present result remains algebraic–operator in nature and does not automatically extend to these applications.

Further work should address extension to general real powers and other rational orders, integral or convolution representations, Laplace and Fourier symbols, initial–boundary value problems, and whether any broader model introduces structure that breaks the similarity equivalence of the normalization family.

## Conclusion

This paper constructs a half-order differentiation operator on a specified graded monomial space directly from the composition requirement. The starting action is

$$D^{\frac{1}{2}} x^{\beta} = c(\beta) x^{\beta - \frac{1}{2}}$$

and requiring two successive half-step actions to reproduce the ordinary derivative gives the recurrence

$$c(\beta) c\left(\beta - \frac{1}{2}\right) = \beta$$

Expanding the recurrence yields two linked coefficient families. For positive integers n,

$$c(n) = c_0 \frac{(2n)!!}{(2n-1)!!}$$

while for the corresponding half-integer values

$$c\left(n - \frac{1}{2}\right) = \frac{(2n-1)!!}{2 c_0 (2n-2)!!}$$

The composition requirement determines the relative coefficients but not the common normalization: reciprocal scale factors in the integer and half-integer sectors cancel after two half-step actions.

The partial Wallis product arising from the double-factorial ratios reveals a structural connection with π but does not determine $c_0$. Compatibility with the Riemann–Liouville monomial formula selects

$$c_0 = \frac{1}{\sqrt{\pi}}$$

was chosen. After this choice, for integer powers we obtain

$$D^{\frac{1}{2}} x^n = \frac{(2n)!!}{(2n-1)!!\sqrt{\pi}} x^{n-\frac{1}{2}}$$

while for the corresponding half-integer powers

$$D^{\frac{1}{2}} x^{n-\frac{1}{2}} = \frac{(2n-1)!!\sqrt{\pi}}{2(2n-2)!!} x^{n-1}$$

On the initial graded space, for monomials and their finite linear combinations,

$$(D^{1/2})^2 f = Df, \quad f \in V_+$$

holds, where $V_+$ is the finite linear span of monomials with non-negative integer and half-integer exponents. For the boundary term generated by applying the first half-step to the constant function, the extension in §1.5 gives

$$D^{\frac{1}{2}} x^{\frac{-1}{2}} = 0$$

which also verifies the composition law for the constant function.

With the chosen normalization, the monomial formula agrees with the Riemann–Liouville rule; for positive integer powers it also agrees with the Caputo rule. For $x^2$, all three calculations yield

$$D_{0+}^{\mathrm{RL},\,1/2} x^2 = D_{0+}^{\mathrm{C},\,1/2} x^2 = D_{\mathrm{alg}}^{1/2} x^2 = \frac{8}{3\sqrt{\pi}} x^{3/2}$$

Agreement on monomials does not imply equivalence on general function spaces. The contribution of the paper is the derivation of the half-order coefficient structure from composition, recurrence, and double factorials within a specified class of linear power-lowering monomial operators. Extension to broader domains, kernel and transform representations, causality, initial conditions, and applied models remains open.